\documentclass{article}
\usepackage[T1]{fontenc}
\fontencoding{T1}
\usepackage{amsfonts}
\usepackage{amssymb}
\usepackage{amsmath}
\newtheorem{theorem}{Theorem}
\newtheorem{lemma}[theorem]{Lemma}

\newcommand{\Norm}{{\rm Norm}}

\newcommand{\Tr}{{\rm Tr}}

\def\qed{\hfill$\square$\\ $\phantom{A}$\\} %\rule{8pt}{8pt}

\title{Symbolic integration in the spirit of Liouville, Abel and Lie}
\author{Waldemar Hebisch}

\begin{document}
\title{Symbolic integration in the spirit of Liouville, Abel and Lie}
\author{
Waldemar Hebisch \\
       {Mathematical Institute}\\
       {Wroc\l{}aw University}\\
       {50-384 Wroc\l{}aw}\\
       {\tt hebisch@mail.math.uni.wroc.pl}
}
\date{}
\maketitle
\begin{abstract}
We provide a Liouville principle for integration in
terms of elliptic integrals.  Our methods are essentially those
of Abel and Liouville changed to modern notation.  We expose
Lie theoretic aspect of Liouville's work.
\end{abstract}
\section{Introduction}

Our task is indefinite integration, that is given
$f$ from some class $A$ we seek $g$ from possibly
different class $B$ such that
$$
f = g'.
$$
As first step in this direction we would like to
delimit form that $g$ can take.  Various results
in this direction are usually called Liouville
principle after famous result of Liouville.

Classically, in computation of field below a hyperbole
it was observed that to integrate rational function one
needs a logarithm.  This lead to class of elementary functions.

In computation of arc length of ellipsis there appear
integral
$$
E(x, m) = \int \frac{(1 - x^2)dx}{\sqrt{(1 - x^2)(1 - mx^2)}}
$$
and question if it can be expressed
% \lquot{}in finite terms\rquot{}
"in finite terms"{}
(in particular question if it is an elementary function).
Abel made first steps and then Liouville proved that
among others elliptic integral $E$ is not elementary.

% To precisely define considered classes of functions we need
% notion of differential field.
Abel and Liouville combined algebraic and analytic arguments,
which causes some doubt about validity of proofs.  Also,
Liouville proof used repeated integration and differentiation
and deeper sense of main part of his argument was not clear.
It turns out that Liouville used differential automorphisms,
integration allowed passing info from generator (a derivation)
to group elements, differentiation went back.  We observed
that this argument can be done directly inside Lie algebra
of derivations.  In effect, our proof is very close to
original Liouville's proof.  Thanks to use of Abel addition
formula for elliptic integrals we can handle elliptic
functions and elliptic integrals.  It is curious that
Liouville missed that: all ingredients were known to him
and our main result is natural extension of that of Liouville.

\section{Differential fields}

Differential field is a field $F$ with a derivation $D$, that
is additive operation which satisfies
Leibniz formula:
$$
D(fg) = D(f)g + fD(g)
$$
We say that $f \in F$ is a constant if and only if $D(f) = 0$.

Remark:  In the sequel we consider only fields of characteristic $0$,
in finite characteristic several crucial results are
no longer valid.

Example: Field ${\mathbb Q}(x, \exp(1/x^2))$.  Function
$f = ((x^2-2)\exp(1/x^2))/x^2$ is an element of this field.

Remark: Given function belongs to many fields, in practice
we prefer small fields.

Example: Field of meromorphic functions in a complex area $U$
with usual derivative is a differential field.

For computational purposes this field is too big, we want
finitely generated fields.
Field of meromorphic functions is in a sense universal example:
every finitely generated differential field is isomorphic to
a subfield of field of meromorphic functions in some area
(Seidenberg \cite{Seid}).

% Example: Formal Laurent or Puiseaux series:
% $$
% \sum_{i = i_0}^\infty a_ix^i
% $$
% where $a_i$ are constants and $Dx = 1$.  This is very big
% field, convenient for some theoretical purposes.

Assumption that we work inside a differential field introduces
some limitations, for example it excludes $|x|$.  In
case of multivalued function sometimes we need to make
choice of branches, for example
$$
\sqrt{x}\sqrt{1-x} - \sqrt{x(1-x)}
$$
is zero or not depending on choice of branches of square root.

In a sense this is necessary limitation, other settings
quickly lead to unsolvable problems (\cite{Rich:U}, \cite{Risch:Prob}).

\subsection{Differential fields, representation}

One way to represent differential field is via generators
and relations.  In characteristic $0$ our field is
an extension of field of rational numbers.  Generators
either are transcendental or defined via minimal
polynomial of an algebraic extension (that is relation).
We define derivative on transcendental generators
in arbitrary way and extend to the whole field
(in characteristic $0$ there exist exactly one extension,
for example see
\cite{Lang:A} case 1 and case 2 following Theorem 5.1).

Field ${\mathbb Q}(x, \exp(1/x^2))$ is given in this way,
there are two generators (both transcendental) $\theta_1 = x$,
$\theta_2 = \exp(1/x^2)$,
$$
D(\theta_1) = 1,$$
$$D(\theta_2) = \frac{-2}{\theta_1^3}\theta_2.$$

Alternative point of view is that a differential field is generated
by solutions of a system of algebraic differential equations.

From this point of view constants are first integrals of our
system of equations, finding them can be a hard problem.

\subsection{Differential fields, geometric model}

Alternatively, we can treat finitely generated differential field
as rational functions on an algebraic manifold $M$.  Derivative
$D$ is a vector field on $M$.

Remark: saying about field of rational functions we automatically
make no distinction between manifolds with the same field
of rational functions (birationaly equivalent).

\subsection{Differential fields, towers}

For us crucial role is played by extensions of transcendental
degree $1$, that is pairs $F \subset K$ of differential fields
such that $K$ is of transcendental degree $1$ over $F$.
Geometrically, $K$ is a field of algebraic functions on
a curve (algebraic manifold of dimension $1$) over $F$.
Assuming that $K$ is algebraic over $F(\theta)$ we get
(variant of) chain rule
$$
Df = D_Ff + (D\theta)\partial_\theta f
$$
where $D_F$ is derivation of $K$ equal to $D$ on $F$ and such
that $D_F(\theta) = 0$ and $\partial_\theta$ is
derivation of $K$ which is zero on $F$ and such that
$\partial_\theta(\theta) = 1$.

We obtain more general fields as towers, that is sequences
$$
F_0 \subset F_1 \subset F_2 \dots \subset F_n = K
$$
where $F_i \subset F_{i+1}$ is an extension of transcendental degree
$\leq 1$.

One can consider more general towers, but here we would like
to consider more specific case.

\subsection{Differential fields, extensions of degree $1$}

Typical examples of extensions of degree $\leq 1$, $F \subset F(\theta)$:
\begin{itemize}
\item algebraic extension
\item extension by a primitive, $D(\theta) = f$ where $f\in F$
\item extension by an exponential, $D(\theta) = D(f)\theta$
where $f\in F$
\item extension by an elliptic function, $D(\theta) = D(f)q$ where
$f\in F$,
$q^2 = \theta^3 -a\theta - b$,  $a$ and $b$ are constants
\item extension by a Lambert W function,
$D(\theta) = D(f)\frac{\theta}{f(\theta+1)}$, where $f\in F$
\end{itemize}
Intuitively in the last three cases we have
$\theta = \exp(f)$, $\theta = P(a, b, f)$, $\theta = W(f)$,
where $P$ is (an alternative version of) Weierstrass $P$ elliptic function
and $W$ is Lambert W function.

In case of elliptic functions corresponding curve is
an elliptic curve, in other cases we have projective line.
However, algebraic extensions can introduce arbitrary curves.

In the sequel we will need two kinds of primitives:
logarithms, where $f = \frac{Dh}{h}$ for $h\in F$ and
elliptic integrals.

\subsection{Differential fields, elliptic integrals}

To define elliptic integrals we assume that there are
elements $p, q \in F$ such that $q^2 = p^3 -ap - b$ where $a$ i $b$ are
constants.

$\theta \in K$ is an elliptic integral of the first kind if
$$
D(\theta) = \frac{Dp}{q}.
$$

$\theta \in K$ is an elliptic integral of the second kind if
$$
D(\theta) = \frac{pDp}{q}.
$$

$\theta \in K$ is an elliptic integral of the third kind if
there exist constant $c$ such that
$$
D(\theta) = \frac{Dp}{(p - c)q}.
$$
Our elliptic integrals are integrals of a differential form
on an elliptic curve in Weierstrass form.

Traditional approach uses Legendre form curve:
$$
y^2 = (1 - x^2)(1 - mx^2).
$$
Those curves give essentially equivalent theories,
however to write a curve in Legendre form we may need
additional algebraic extensions.  For our purposes we need
control over algebraic extensions so Weierstrass form
is preferable.  Also, integrals on Weierstrass curve fit
well with Weierstrass elliptic functions.  Legendre curve
naturally leads to  Jacobi elliptic functions which for
us are less convenient.

Elliptic integrals in Legendre form are frequently
transformed to trigonometric form.  For us this is
complication which needlessly introduces transcendental
functions.

Crucial for us Abel lemma is proved on curve in Legendre
form, thanks to equivalence we will get it also for
Weierstrass elliptic integrals.

\subsection{Differential fields, Lie closed extensions}

All extensions of degree $1$ that we considered are
Lie closed, that is there exists nonzero derivation $X$
on $K$ such that $X$ commutes with $D$
and $X$ is zero on $F$.
More generally extension $F \subset K$ of transcendental degree $n$
is Lie closed
if there are $n$ linearly independent derivations ${X_k}$ on
$K$ which commute with $D$ and are zero on $F$.

Explicitly:
\begin{itemize}
\item For $\theta$ which is a primitive $X = \partial_\theta$
\item For $\theta$ which is an exponential
$X = \theta\partial_\theta$
\item For $\theta$ which is an elliptic function, that is
$D(\theta) = D(f)q$ where $f\in F$, $q^2 = \theta^3 -a\theta - b$,
$a$ and $b$ are constants we take $X = q\partial_\theta$
\item For $\theta$ which is a Lambert W
$$
X = \frac{\theta}{\theta + 1}\partial_\theta
$$
\end{itemize}

Note that $DX - XD$ is a derivation, so it is enough to
check commutation on generators, that is elements of $F$ and $\theta$.
On  $F$ our $X$ is zero, $D$ preserves $F$, so both products are
zero, so $D$ and $X$ commute on $F$.

For $\theta$ which is a primitive we take $X = \partial_\theta$
so we have $\partial_\theta\theta = 1$, so $DX\theta = 0$.
$D\theta \in F$ so $XD\theta = 0$ and again $D$ and $X$ commute.

More generally, for all our $\theta$ we have $D\theta = gh(\theta)$
where $g \in F$
and $h(\theta)$ depends only on $\theta$ (in the elliptic case
$h$ is an algebraic function, in other cases it is rational function).
Now our $X$ is
$$
X = h(\theta)\partial_\theta
$$
and we have
$$
XD\theta = Xgh(\theta) = gXh(\theta) =
gh(\theta)\partial_\theta h(\theta),
$$
$$
DX\theta = Dh(\theta) = gh(\theta)\partial_\theta h(\theta)
$$
so $D$ and $X$ commute.

In three cases $X$ generates one parameter group of
(differential) automorphisms:
\begin{itemize}
\item When $\theta$ is a primitive, then map
$\theta \mapsto \theta + c$ where $c$ is a constant extends
to an automorphism of $K$.  In other words translation by
a constant gives automorphism of differential field.
\item When $\theta$ is an exponential, then mapping
$\theta \mapsto c\theta$ where $c$ is nonzero constant extends
to an automorphism.  So group is multiplicative group of
nonzero constants.
\item On elliptic curve we have multiplication, multiplication
by points with constant coordinates gives automorphisms.
\end{itemize}
In other words extensions above are differential Galois
(Kolchin proves that there are no other differential Galois
extensions of degree $1$).

The case of Lambert W is different: on algebraic level
there are no group.

\subsection{Differential fields, indefinite integral}

We sat that $g \in F$ is indefinite integral (or a primitive)
of $f \in F$ when
$$
f = D(g).
$$
For example
$$
\int ((x^2-2)\exp(1/x^2))/x^2 = x\exp(1/x^2).
$$
As element of differential field indefinite integral is
uniquely determined up to additive constant
(but element of a differential field may be represented
by multiple expressions, so we can get different formulas).

\subsection{Elementary extensions}
Differential field $K$ is an elementary extension of $F$
if there exists tower
$$
F = F_0 \subset F_2 \subset \dots \subset F_n = K
$$
where each of $F_i \subset F_{i+1}$ is an algebraic extension,
extension by a logarithm or extension by an exponential.

We say that a function $f$ is elementary over $F$ if it
is an element of an elementary extension of $F$.

We say that a function $f$ is elementary if it is elementary
over field of rational functions over constants.

Example:  Let $f(x) = \sqrt{\exp(x + \log(x))}$.
$f \in K = {\mathbb Q}(x)(\theta_1, \theta_2, \theta_3)$ where
$\theta_1 = \log(x)$, $\theta_2 = \exp(x + \log(x))$,
$\theta_3 = \sqrt{\exp(x + \log(x))}$.  So $f$ is an
elementary function.

We can write function from previous example as
$f(x) = \sqrt{x\exp(x)}$,
and use $K = {\mathbb Q}(x)(\theta_1, \theta_2)$ for
$\theta_1 = \exp(x)$, $\theta_2 = \sqrt{x\exp(x)}$.
However, given an expression for $f(x)$ we can build
an elementary extension in a natural way: each logarithm
and exponential appearing in $f$ and each irrational algebraic
subexpression of $f$ (like a root) is associated with a
generator of an extension.  Different expressions for $f$
may lead to different elementary extensions.

We normally assume that extensions can not be written in
simpler form.  For example we will treat $\theta$ as
an exponential or a logarithm only when it is not algebraic
over smaller field.  For example we treat
$f(x) = \exp(\frac{\log(x)}{2}) = \sqrt{x}$
not as an exponential but as a solution to algebraic
equation $f^2(x) = x$.  Similarly we do not treat
$x = \log(\exp(x))$ as a logarithm.

\subsection{Elliptic-Lambert extensions}

We will consider also wider class of extensions
"{}elliptic-Lambert"{} extension where we allow
towers in which may appear extensions by elliptic functions,
elliptic integrals and Lambert W function.

\subsection{Commutator formula}

\begin{lemma}\label{der-comm}
If $X, Y$ are derivations on $K$, $p\in K$, $p$ and $\psi(p)$ are
algebraically dependent over constants, then
$$
X((Yp)\psi(p)) - Y((Xp)\psi(p) = ([X, Y]p)\psi(p)
$$
\end{lemma}

Proof: Since derivations extend uniquely to algebraic extensions
we have $X\psi(p) = (Xp)\psi'(p)$ and $Y\psi(p) = (Yp)\psi'(p)$
where ${}'$ denotes unique derivation on $C(p, \psi(p))$ such
that $p' = 1$ and $C$ is constant field.  Now,
$$
X((Yp)\psi(p)) = (XYp)\psi(p) + (Yp)(Xp)\psi'(p),
$$
$$
Y((Xp)\psi(p) = (YXp)\psi(p) + (Xp)(Yp)\psi'(p).
$$
Subtracting the above give the result.
\qed.

Remark: While different proofs of Lioville's theorem at first
may look quite different, crucial step of known proofs either
explicitly uses something like Lemma \ref{der-comm} (for
example Lemma inside proof of Theorem 1 in \cite{BK}) or
depend on calculations which work only because the Lemma is valid
(as is the case with original Liouville's proof).

\section{Abel formula}

We will need Abel addition formula for elliptic integrals
(in Legendre form).  Abel work \cite{Abel:prec} can be
easily modified to modern standards.  However, below we present
somewhat different argument.

 We consider integrals of differential forms
on a curve $C$ in Legendre form
$y^2 = (1 - x^2)(1 - mx^2)$.
Then
$$
\Pi'(x) = \frac{1}{(1 - nx^2)y}
$$
Let $(x_1, y_1)$ and $(x_2, y_2)$ be points on $C$
and $(x_3, y_3)$ be their sum.  Abel gave formula
$$
\Pi(x_1) + \Pi(x_2) = C + \Pi(x_3) -
\frac{a}{2\Delta(a)}\log\left(
\frac{a_0a + a^3 + x_1x_2x_3\Delta(a)}{a_0a + a^3 - x_1x_2x_3\Delta(a)}
\right)
$$
where
$$
n = \frac{1}{a^2},
$$
$$\Delta(a) = \sqrt{(1 - a^2)(1 - ma^2)},$$
$$a_0 = \frac{x_2^3y_1 - x_1^3y_2}{x_1y_2 - x_2y_1}.$$

In differential version
$$
\frac{dx_1}{(1 - x_1^2/a^2)y_1} + \frac{dx_2}{(1 - x_2^2/a^2)y_2}
= \frac{dx_3}{(1 - x_3^2/a^2)y_3} - \frac{a df}{2\Delta(a) f}
$$
where
$$
f = \frac{a_0a + a^3 + x_1x_2x_3\Delta(a)}{a_0a + a^3 - x_1x_2x_3\Delta(a)}.
$$
In this version formula can be checked by direct calculation
(it can be done on a computer), we skip details.

We will also need integrals of the first kind $F$ and
of second kind $E$:
$$
F(x)' = \frac{1}{y},
$$
$$
E(X)' = \frac{1-mx^2}{y}.
$$

\begin{lemma}\label{el-add}(Abel)
$$
\sum_{k=1}^l F(x_i)' = F(y)',
$$
$$
\sum_{k=1}^l E(x_i)' = E(y)' + g',
$$
$$
\sum_{k=1}^l \Pi(x_i)' = \Pi(y)' + f',
$$
where $y$ is sum of points on curve $C$, $g$ is a rational
function and $f$ is a sum of
logarithms.
\end{lemma}

Proof: Formula for $\Pi$ follows by induction from the formula
above.  Formulas for $F$ and $E$ are obtained in similar
way: direct calculation for $l=2$ and induction.  Explicit
formula for $g$ when $l = 2$ is:
$$
g = m\frac{-x_1 x_2^3+x_1^3 x_2}{x_1 y_2 -x_2 y_1} 
$$
\qed

Remark: Abel obtained formula for $F$ taking limit of formula for $\Pi$
when $n$ goes to $0$ (and similarly for $E$).  This can
be justified in purely algebraic way, but in computer era
direct calculation is simpler.

\section{Liouville-Ostrowski theorem}

Liouville-Ostrowski theorem:
\begin{theorem}
If
$F$ is a differential field,
$f\in F$ and $f$ has a primitive in an elementary extension $K$
of $F$, then there exists extension $\bar F$ of $F$ by algebraic
constants and functions
$v_i \in \bar F$ and constants $c_1, \dots, c_l \in \bar F$
such that
$$
f = D(v_0) + \sum_{i=1}^l c_i\frac{D(v_i)}{v_i} =
D(v_0) + \sum_{i=1}^l c_iD(\log(v_i))
$$
\end{theorem}
For modern proof see \cite{Ros:L} (and \cite{Risch:Prob}
to get strong result about constants).

Liouville-Ostrowski theorem says that we can find all parts of integral
of $f$ already in $F$ extended by algebraic constants.  This is
crucial property for symbolic integration algorithms.

Using his theorem Liouville proved that $\int \frac{\exp(x)}{x}dx$,
$\int \exp(x^2)dx$ and elliptic integrals of the first and
second kind are not elementary.

% Example:
% $$
% \int \frac{dx}{\sqrt{x^3 - ax - b}}
% $$
% is not elementary.  Differential form above is locally analytic,
% without any singularities.  Each pole of $v_0$ gives
% pole of order bigger than $1$ in the sum which can not cancel
% with other terms because poles of logarithmic terms are are of
% order $1$.  So $v_0$ has no singularities, so it is
% a constant.  Similarly, we can rewrite logarithmic terms
% so that there is no cancellation between poles, so $v_i$
% have no singularities, so they are constants.

% Example:
% $$
% \int \frac{\exp(x)}{x}dx
% $$
% and
% $$
% \int \exp(x^2)dx
% $$
% are not elementary.  In both cases considering poles in
% w $M(\theta)$ where $M={\mathbb Q}(x)$ and $\theta = \exp(x)$
% (respectively $\exp(x^2)$) shows that integral have form
% $$
% R(x)\theta
% $$
% Now somewhat tedious but straightforward reasoning in ${\mathbb Q}(x)$
% shows that there are no solutions.

\subsection{Generalization}
\begin{theorem}\label{liu-ab}
If
$F$ is a differential field, $f\in F$ and $f$ has a primitive
in an elliptic-Lambert extension $K$ of $F$, then
there exist an extension $\bar F$ of $F$ by algebraic constants,
functions
$v_i \in \bar F$ and constants $c_1, \dots, c_l \in \bar F$, such that
$$
f = D(v_0) + \sum_{i=1}^l c_i\phi(Dv_i, v_i)
$$
where $\phi(Dv_i, v_i)$ is of one of forms below
\begin{itemize}
\item derivative of a logarithm, that is
$\phi(Dv_i, v_i) = \frac{D(v_i)}{v_i}$
\item $\phi(Dv_i, v_i) = \frac{Dv_i}{q_i}$
\item $\phi(Dv_i, v_i) = \frac{v_iDv_i}{q_i}$
\item $\phi(Dv_i, v_i) = \frac{Dv_i}{(v_i - c_i)q_i}$
\end{itemize}
where $q_i^2 = v_i^3 - a_iv_i - b_i$, $q_i \in \bar F$, $a_i, b_i, c_i$
are constants in $\bar F$.
\end{theorem}

In other word, we can find all ingredients of the integral already in $F$
extended by algebraic constants.

Remark: Abel proved equivalent result in case when $F$ is
algebraic over $C(x)$ (where $C$ is constant field) and
$K$ is algebraic over $F$.

Proof:
Proof is via induction over tower.  Namely let $L$ be
elliptic-Lambert extension in which integral exists.
We can write
$$
\bar F = F_0 \subset F_1 \subset \dots \subset F_n = L
$$
where each $F_j \subset F_{j+1}$ is an extension
of degree $ \leq 1$ of form given earlier.
By assumption in $L$ we can write integral in the form
given in the theorem.  So it remains to prove that
given expression as above with all parts in $F_{j+1}$
we can transform it into expression with all parts
in $F_j$ extended by constants.  We will do this in a few steps.

Step 1.  We will prove that:
$$
X\frac{Dv}{v} = D\frac{Xv}{v}
$$
where $X$ is given earlier derivation commuting with $D$.
Also
$$
X\frac{Dp}{(p - c)q} = D\frac{Xp}{(p - c)q}
$$
and similarly for remaining $\phi(Dv_i, v_i)$ terms.

Namely, all our $\phi(Dv_i, v_i)$ are of form
$(Dv_i)\psi(v_i)$ where $\psi$ is an algebraic function.
Since $[X, D] = 0$ equality
$$
X\phi(Dp, p) = D\phi(Xp, p).
$$
follows from Lemma \ref{der-comm}.

Step 2. Now we compute
$$
0 = Xf = XD(v_0) + \sum_{i=1}^l c_iX\phi(Dv_i, v_i)
$$
$$
= DXv_0 + \sum_{i=1}^l c_iD\phi(Xv_i, v_i)
$$
$$
= D(Xv_0 + \sum_{i=1}^l c_i\phi(Xv_i, v_i))
$$
So
$$
c = Xv_0 + \sum_{i=1}^l c_i\phi(Xv_i, v_i)
$$
is a constant.

Step 3. Now we use chain rule:
$$
D = D_{F_j} + (D\theta)\partial_\theta
$$
If $\theta$ is a primitive, then the above can be written as
$$
D = D_{F_j} + (D\theta)X.
$$

Now
$$
D(v_0) + \sum_{i=1}^l c_i\phi(Dv_i, v_i) =
D_{F_j}(v_0) + \sum_{i=1}^l c_i\phi(D_{F_j}v_i, v_i) +
$$
$$
(D\theta)(Xv_0 + \sum_{i=1}^l c_i\phi(Xv_i, v_i))
$$
$$
= D_{F_j}(v_0) + \sum_{i=1}^l c_i\phi(D_{F_j}v_i, v_i) + cD\theta
$$
But $D\theta = \phi(Dp, p)$ is in the form required by the theorem, so
we can add it as an additional term in the sum.
$D_{F_j}(\theta) = 0$, so we have expression of required form in
$F_j$ extended by constants.

Step 3'.  In Lambert W case we have $D\theta = \frac{D(v)}{v}X\theta$
so $D = D_{F_j} + \frac{D(v)}{v}X$ and $cD\theta$ term can be replaced
by logarithmic term $\frac{D(v)}{v}$.

Step 3''.  In other cases (extension by exponential or
elliptic function) $D\theta = D(v)X\theta$
and proceeding as before we add $cD(v)$ to $v_0$.

It remains to consider algebraic extensions.  We do this using
methods of Galois theory (or more precisely method of Abel).
We use trace and norm maps.  For an algebraic extension
$F \subset E$ we define
$$
\Tr(y) = \sum_{k=1}^{m}\iota_k(y),
$$
$$
\Norm(y) = \prod_{k=1}^{m}\iota_k(y)
$$
where $\iota_k$ goes over all embeddings of $E$ over $F$
into algebraic closure of $F$.

We have
$$
\frac{D\Norm(y)}{\Norm(y)} = \Tr(\frac{Dy}{y}).
$$
Namely
$$\Norm(y) = \prod_{k=1}^{m} \iota_k(y).$$
Since derivative extends uniquely onto algebraic extension
we have
$D\iota_k(y) = \iota_k(Dy)$.
By the logarithmic derivative formula we have
$$
\frac{D\Norm(y)}{\Norm(y)} = \sum \frac{D\iota_k(y)}{\iota_k(y)}
= \sum \frac{\iota_k(Dy)}{\iota_k(y)} =  \sum \iota_k(\frac{D(y)}{y})
= \Tr(\frac{Dy}{y}).
$$

From Abel formula (Lemma \ref{el-add})
$$
\Tr(D(\Pi(p)) = D\Pi(\tilde p) + Df
$$
where $\tilde p$ is sum of images on curve and $f$ is a sum of
logarithms.  Namely, put $p_k = \iota_k(p)$.  We have
$$
D\Pi(p_k) = D\Pi(\iota_k(p)) = \iota_k(D\Pi(p))
$$
so applying Lemma \ref{el-add} we get
$$
\Tr(D(\Pi(p)) = \sum \iota_k(D(\Pi(p))) = \sum D\Pi(\iota_k p))
$$
$$
= D\Pi(\tilde p) + Df.
$$
Note that Lemma \ref{el-add} is written in terms of
derivative $D(f) = f'$ but is really equality of differential forms
so one can substitute an arbitrary derivative and equality
remains valid.

  By Galois (Abel) theory $\tilde p$ is in $E$.
Similar result holds for elliptic integrals of the first and
second kind.  Passing between  Legendre and Weierstrass theory
we get similar result for integrals in Weierstrass form.

Now, when $F_j \subset F_{j+1}$ is an algebraic extension
and in $F_{j+1}$ we have
$$
f = D(v_0) + \sum_{i=1}^l c_i\phi(Dv_i, v_i)
$$
and $f \in F_j$, then applying $\Tr$ to both sides we get similar
formula with terms in $F_j$ extended by constants, which
ends the proof of algebraic case.

Our proof can introduce transcendental constants.  We should
prove that it is enough to use algebraic constants.  This
can be done in various ways, for example using Hilbert theorem
about zeros (like \cite{Risch:Prob}) or adapting model-theoretic
proof of Hilbert theorem.
We will skip details.
\qed

\section{Further results and remarks}

We can strengthen our main theorem by allowing more complicated $K$.
Namely, we can allow towers containing Lie closed extensions
such that Lie algebra spanned by $X_k$ is spanned by commutators.
Put
$$w_{X_k} = X_kv_0 + \sum_{i=1}^l c_i\phi(X_kv_i, v_i).$$
Proceeding like in step 2 of proof of Theorem we get:
$$
X_kw_{X_j} - X_jw_{X_k} =
[X_k, X_j]v_0 + \sum_{i=1}^l c_i\phi([X_k, X_j]v_i, v_i).
$$
From step 2 we know that $w_{X_k}$ is a constant, so
$$
X_kw_{X_j} - X_jw_{X_k} = 0
$$
hence
$$
[X_k, X_j]v_0 + \sum_{i=1}^l c_i\phi([X_k, X_j]v_i, v_i) = 0
$$
We assume that Lie algebra spanned by $X_k$ is spanned by
commutators, so also
$$w_{X_k} = 0.$$
But having this we can proceed like in step 3 using formula
$$
D = D_{F_j} + \sum a_kX_k
$$
where $a_k \in K$ are such that $D\theta_i = \sum a_k X_k\theta_i$
where ${\theta_i}$ is transcendence basis of $F_{j+1}$ over $F_j$.
Terms involving $X_k$ will vanish so
$$
f = Dv_0 + \sum_{i=1}^l c_i\phi(Dv_i, v_i) =
    D_{F_j} + \sum_{i=1}^l c_i\phi(D_{F_j}, v_i).
$$

It is well known (see for example \cite{Nish:Lie}) that
Picard-Vessoit extensions are Lie closed and that semisimple
Lie algebra is generated by commutators.  Many classical
special functions are solutions of linear ordinary differential
equations, so are elements of Picard-Vessoit extensions.
Generically corresponding Lie algebra is semisimple.  Let
us note that there are a lot of integrals which can be
expressed in terms of hypergeometric functions but not in
terms of elementary functions.  Our result proves that
all those hypergeometric functions must correspond to
degenerate cases of hypergeometric equation.

Theorem \ref{liu-ab} is about functions integrable
in terms of elliptic integrals.  However, if $f$ is integrable
in some extension $L$ of $F$ which is last term of tower
$$
F = F_0 \subset F_1 \subset \dots \subset F_n = L
$$
where each of $F_i \subset F_{i+1}$ is either elementary
or extension by Lambert W
or extension by elliptic function or Picard-Vessoit extensions
with semisimple Lie algebra, then $f$ has elementary integral.
Namely, in our proof elliptic integrals appeared only because
extension by elliptic integral was part of the tower.

Result about elliptic integrals of first and second kind can
be proved via main theorem in \cite{BK}, however to
handle elliptic integrals of third kind we need Abel
formula.

\end{document}